\documentclass[a4paper,11pt]{article}

\usepackage{microtype}

\usepackage[T1]{fontenc}
\usepackage[utf8]{inputenc}

\usepackage{amsfonts}
\usepackage{amsmath}
\usepackage{amssymb}
\usepackage{mathtools}
\usepackage[binary-units=true]{siunitx}
\usepackage[boldvectors,braket]{physymb}
\usepackage{bm}

\DeclareMathOperator{\Perm}{Perm}

\usepackage{caption}
\usepackage{subcaption}
\usepackage{floatrow}
\usepackage{algpseudocode}
\DeclareNewFloatType{algorithm}{name=Algorithm}
\usepackage{txfonts}

\usepackage{multicol}
\usepackage{booktabs}

\makeatletter
\newcommand\tabfill[1]{%
\dimen@\linewidth%
\advance\dimen@\@totalleftmargin%
\advance\dimen@-\dimen\@curtab%
\parbox[t]\dimen@{#1\ifhmode\strut\fi}%
}

\usepackage{color}
\definecolor{thered}{rgb}{0.65,0.04,0.07}
\definecolor{thegreen}{rgb}{0.06,0.44,0.08}
\definecolor{theblue}{rgb}{0.02,0.2,0.68}

\usepackage[numbers,sort,compress]{natbib}
\usepackage[unicode=true,
            bookmarks=true,
            bookmarksnumbered=false,
            bookmarksopen=false,
            breaklinks=true,
            pdfborder={0 0 1},
            backref=false,
            colorlinks=true,
            linkcolor=black,
            urlcolor=theblue,
            citecolor=theblue,
            hyperfootnotes=false]
           {hyperref}
\hypersetup{pdftitle={},
            pdfauthor={F. D. Witherden, P. E. Vincent}}

\begin{document}

\title{On the Identification of Symmetric Quadrature Rules for Finite
  Element Methods}

\author{F. D. Witherden\footnote{Corresponding author; e-mail
    freddie.witherden08@imperial.ac.uk.}, P. E. Vincent\\\\
  \textit{\small Department of Aeronautics, Imperial College London, SW7
    2AZ}}
\maketitle

\begin{abstract}
  In this paper we describe a methodology for the identification of
symmetric quadrature rules inside of quadrilaterals, triangles,
tetrahedra, prisms, pyramids, and hexahedra. The methodology is free from
manual intervention and is capable of identifying an ensemble of rules
with a given strength and a given number of points. We also present
\emph{polyquad} which is an implementation of our methodology. Using
polyquad we proceed to derive a complete set of symmetric rules on the
aforementioned domains. All rules possess purely positive weights and have
all points inside the domain. Many of the rules appear to be new, and an
improvement over those tabulated in the literature.
\end{abstract}

\newpage

\section{Introduction}

When using the finite element method to solve a system of partial
differential equations it is often necessary to evaluate surface and
volume integrals inside of a standardised domain $\vec{\Omega}$
\cite{hesthaven2008nodal, solin2003higher, karniadakis2013spectral}.  A
popular numerical integration technique is that of Gaussian quadrature
in which
\begin{equation} \label{eq:quad}
  \int_{\vec{\Omega}} f(\vec{x})\, \mathrm{d} \vec{x} \approx \sum_i^{N_p}
  \omega_i f(\vec{x}_i),
\end{equation}
where $f(\vec{x})$ is the function to be integrated, $\{\vec{x}_i\}$ are
a set of $N_p$ points, and $\{\omega_i\}$ the set of associated weights.
The points and weights are said to define a \emph{quadrature rule}.  A
rule is said to be of strength $\phi$ if it is capable of exactly
integrating any polynomial of maximal degree $\phi$ over $\vec{\Omega}$.
A degree $\phi$ polynomial $p(\vec{x})$ with $\vec{x} \in \vec{\Omega}$
can be expressed as a linear combination of basis polynomials
\begin{equation}
  p(\vec{x}) = \sum_i^{|\mathcal{P}^\phi|} \alpha_i
  \mathcal{P}^{\phi}_i(\vec{x}),  \qquad
  \alpha_i = \int_{\vec{\Omega}} p(\vec{x}) \mathcal{P}^{\phi}_i(\vec{x})\, \mathrm{d} \vec{x},
\end{equation}
where $\mathcal{P}^\phi$ is the set of basis polynomials of degree $\leq
\phi$.  From the linearity of integration it therefore follows that a
strength $\phi$ quadrature rule is one which can exactly integrate the
basis.  Taking $f \in \mathcal{P}^\phi$ the task of obtaining an $N_p$
point quadrature rule of strength $\phi$ is hence reduced to finding a
solution to a system of $|\mathcal{P}^\phi|$ nonlinear equations.  This
system can be seen to possess $(N_D + 1)N_p$ degrees of freedom where
$N_D \geq 2$ corresponds to the number of spatial dimensions.

In the case of $N_p \lesssim 10$ the above system can often be solved
analytically using a computer algebra package.  However, beyond this it
is usually necessary to solve the above system---or a simplification
thereof---numerically.  Much of the research into multidimensional
quadrature over the past five decades has been directed towards the
development of such numerical methods. The traditional objective when
constructing quadrature rules is to obtain a rule of strength $\phi$
inside of a domain $\vec{\Omega}$ using the fewest number of points.  To
this end efficient quadrature rules have been derived for a variety of
domains: triangles \cite{lyness1975moderate, dunavant1985high,
  lyness1994survey, savage1996quadrature, wandzurat2003symmetric,
  zhang2009set, taylor2005several, xiao2010numerical,
  witherden2013analysis, williams2014symmetric}, quadrilaterals
\cite{dunavant1985economical, cools1988another, xiao2010numerical},
tetrahedra \cite{savage1996quadrature, zhang2009set, shunn2012symmetric,
  keast1986moderate}, prisms \cite{kubatko2013pri}, pyramids
\cite{kubatko2013pyr}, and hexahedra \cite{stroud1971approximate,
  dunavant1986efficient, cools2001rotation, xiao2010numerical}.  For
finite element applications it is desirable that (i) points are arranged
symmetrically inside of the domain, (ii) all of the points are strictly
inside of the domain, and (iii) all of the weights are positive.  The
consideration given to these criterion in the literature cited above
depends strongly on the intended field of application---not all rules
are derived with finite element method in mind.

Much of the existing literature is predicated on the assumption that the
integrand sits in the space of $\mathcal{P}^\phi$.  Under this
assumption there is little, other than the criteria listed above, to
distinguish two $N_p$ rules of strength $\phi$; both can be expected to
compute the integral exactly with the same number of functional
evaluations.  It is therefore common practice to terminate the rule
discovery process as soon as a rule is found.  However, there are cases
when either the integrand is inherently non-polynomial in nature, e.g.
the quotient of two polynomials, or of an high degree, e.g. a polynomial
raised to a high power.  In these circumstances the above assumption no
longer holds and it is necessary to consider the truncation term
associated with each rule.  Hence, within this context it is no longer
clear that the traditional objective of minimising the number of points
required to obtain a rule of given strength is suitable: it is possible
that the addition of an extra point will permit the integration of
several of the basis functions of degree $\phi + 1$.

Over the past five or so years there has also been an increased interest
in numerical schemes where the same set of points are used for both
integration and interpolation.  One example of such a scheme is the flux
reconstruction (FR) approach introduced by Huynh \cite{huynh2007flux}.
In the FR approach there is a need for quadrature rules that (i) are
symmetric, (ii) remain strictly inside of the domain, (iii) have a
prescribed number of points, and (iv) are associated with a well
conditioned nodal basis for polynomial interpolation.  These last two
requirements exclude many of the points tabulated in the literature.
Consequently, there is a need for \emph{bespoke} or \emph{designer}
quadrature rules with non-standard properties.

This paper describes a methodology for the derivation of symmetric
quadrature rules inside of a variety of computational domains.  The
method accepts both the number of points and the desired quadrature
strength as free parameters and---if successful---yields an ensemble of
rules.  Traits, such as the positivity of the weights, can then be
assessed and rules binned according to their suitability for various
applications.  The remainder of this paper is structured as follows.  In
\autoref{sec:shapes} we introduce the six reference domains and
enumerate their symmetries.  Our methodology is presented in
\autoref{sec:meth}.  Based on the approach of Witherden and Vincent
\cite{witherden2013analysis} the methodology requires no manual
intervention and avoids issues relating to ill-conditioning.  In
\autoref{sec:impl} we proceed to describe our open-source
implementation, \emph{polyquad}.  Using polyquad a variety of
truncation-optimised rules, many of which appear to improve over those
tabulated in the literature, are obtained and presented in
\autoref{sec:rules}.  Finally, conclusions are drawn in
\autoref{sec:conclusions}.

\section{Bases, Symmetries, and Domains}
\label{sec:shapes}

\subsection{Basis polynomials}

The defining property of a quadrature rule for a domain $\vec{\Omega}$
is its ability to exactly integrate the set of basis polynomials,
$\mathcal{P}^\phi$.  This set has an infinite number of representations
the simplest of which being the monomials.  In two dimensions we can
express the monomials as
\begin{equation}
  \mathcal{P}^\phi = \bigl\{x^iy^j \mid 0 \leq i \leq \phi,\; 0 \leq j
  \leq \phi - i \bigr\},
\end{equation}
where $\phi$ is the maximal degree.  Unfortunately, at higher degrees
the monomials become extremely sensitive to small perturbations in the
inputs.  This gives rise to polynomial systems which are poorly
conditioned and hence difficult to solve numerically \cite{zhang2009set,
  shunn2012symmetric}.  A solution to this is to switch to an
\emph{orthonormal basis set} defined in two dimensions as
\begin{equation}
  \mathcal{P}^\phi = \bigl\{\psi_{ij}(\vec{x}) \mid 0 \leq i \leq \phi,\; 0 \leq j
  \leq \phi - i \bigr\},
\end{equation}
where $\vec{x} = (x,y)^T$ and $\psi_{ij}(\vec{x})$ is satisfies $\forall
\mu, \nu$

\begin{equation}
  \int_{\vec{\Omega}} \psi_{ij}(\vec{x}) \psi_{\mu\nu}(\vec{x})\,
  \mathrm{d}\vec{x} = \delta_{i\mu}\delta_{j\nu},
\end{equation}
where $\delta_{i\mu}$ is the Kronecker delta.  In addition to being
exceptionally well conditioned orthonormal polynomial bases have other
useful properties.  Taking the constant mode
of the basis to be $\psi_{00}(\vec{x}) = 1/c$ we see that
\begin{equation} \label{eq:obint}
  \int_{\vec{\Omega}} \psi_{ij}(\vec{x}) \, \mathrm{d}\vec{x} =
  c \int_{\vec{\Omega}} \psi_{00}(\vec{x}) \psi_{ij}(\vec{x})
  \, \mathrm{d}\vec{x} = c\delta_{i0}\delta_{j0},
\end{equation}
from which we conclude that all non-constant modes of the basis
integrate up to zero.  Following Witherden and Vincent
\cite{witherden2013analysis} we will use this property to define the
truncation error associated with an $N_p$ point rule
\begin{equation} \label{eq:trunc}
  \xi^2(\phi) = \sum_{i,j}
  \bigg\{\sum_{k}^{N_{p}}\omega_k \psi_{ij}(\vec{x}_k) -
  c\delta_{i0}\delta_{j0}\bigg\}^2,
\end{equation}
This definition is convenient as it is free from both integrals and
normalisation factors.  The task of constructing an $N_p$ point
quadrature rule of strength $\phi$ is synonymous with finding a set of
points and weights that minimise $\xi(\phi)$.

Although the above discussion has been presented primarily in two
dimensions all of the ideas and relations carry over into three
dimensions.

\subsection{Symmetry orbits}

A symmetric arrangement of $N_p$ points inside of a reference domain can
be decomposed into a linear combination of \emph{symmetry orbits}.  This
concept is best elucidated with an example.  Consider a line segment
defined by $[-1,1]$.  The segment possesses two symmetries: an identity
transformation and a reflection about the origin.  For an arrangement of
distinct points to be symmetric it follows that if there is a point at
$\alpha$ where $0 < \alpha \leq 1$ there must also be a point at
$-\alpha$.  We can codify this by writing $\mathcal{S}_2(\alpha) =
\pm\alpha$ with $\abs{\mathcal{S}_2} = 2$.  The function $S_2$ is an
example of a \emph{symmetry orbit} that takes a single \emph{orbital
  parameter}, $\alpha$, and generates two distinct points.  In the limit
of $\alpha \rightarrow 0$ the two points become degenerate.  We handle
this degeneracy by introducing a second orbit, $\mathcal{S}_1 = 0$, with
$\abs{\mathcal{S}_1} = 1$.  Having identified the symmetries we may now
decompose a symmetric arrangement of points as
\[
  N_p = n_1\abs{\mathcal{S}_1} + n_2\abs{\mathcal{S}_2} = n_1 + 2n_2,
\]
where $n_1 \in \{0, 1\}$ and $n_2 \geq 0$ with the constraint on $n_1$
being necessary to ensure uniqueness.  This is a constrained linear
Diophantine equation; albeit one that is trivially solvable and admits
only a single solution.  As a concrete example we take $N_p = 11$.
Solving the above equation we find $n_1 = 1$ and $n_2 = 5$.  The $n_1$
orbit does not take any arguments and so does not contribute any degrees
of freedom.  Each $n_2$ orbit takes a single parameter, $\alpha$, and so
contributes one degree of freedom for a grand total of five.  This is
less than half that associated with the asymmetrical case.  Hence, by
parameterising the problem in terms of symmetry orbits it is possible to
simultaneously reduce the number of degrees of freedom while
guaranteeing a symmetric distribution of points.

Symmetries also serve to reduce the number of basis polynomials that
must be considered when computing $\xi(\phi)$.  Consider the following
two monomials
\[
 p_1(x, y) = x^iy^j, \qquad \text{and} \qquad p_2(x, y) = x^jy^i,
\]
defined inside of a square domain with vertices $(-1,-1)$ and $(1,1)$.
We note that $p_1(x,y) = p_2(y,x)$.  As this is a symmetry which is
expressed by the domain it is clear that any symmetric quadrature rule
capable of integrating $p_1$ is also capable of integrating $p_2$.
Further, the index $i$ is odd we have $p_1(x,y) = -p_1(-x,y)$.
Similarly, when $j$ is odd we have $p_1(x,y) = -p_1(x, -y)$.  In both
cases it follows that the integral of $p_1$ is zero over the domain.
More importantly, it also follows that \emph{any} set of symmetric
points are also capable of obtaining this result.  This is due to terms
on the right hand side of \autoref{eq:quad} pairing up and cancelling
out.  A consequence of this is that not all of equations in the system
specified by \autoref{eq:quad} are independent.  Having identified such
polynomials for a given domain it is legitimate to exclude them from our
definition of $\xi(\phi)$.  Although this exclusion does change the
value of $\xi(\phi)$ in the case of a non-zero truncation error the
effect is not significant.  We shall denote the set of basis polynomials
which \emph{are} included as the \emph{objective basis}, and denote this
by $\mathcal{\tilde{P}}^\phi$.

\subsection{Reference domains}

In the paragraphs which follow we will take
$\hat{P}^{(\alpha,\beta)}_i(x)$ to refer to a \emph{normalised} Jacobi
polynomial as specified in \S 18.3 of \cite{olver2010nist}.  In two
dimensions we take the coordinate axes to be $\vec{x} = (x,y)$ and
$\vec{x} = (x,y,z)$ in three dimensions.

\paragraph{Triangle.}

\begin{figure}
  \centering
  \begin{subfigure}[b]{.45\linewidth}
    \centering
    \includegraphics{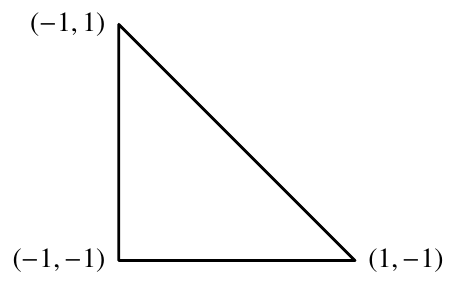}
    \caption{Triangle.}
  \end{subfigure}
  \begin{subfigure}[b]{.45\linewidth}
    \centering
    \includegraphics{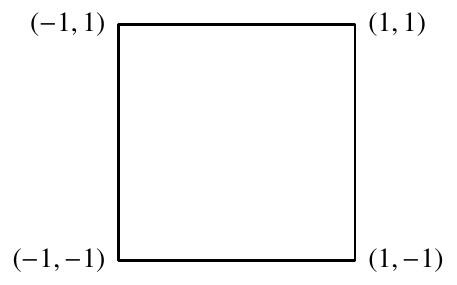}
    \caption{Quadrilateral.}
  \end{subfigure}
  \caption{\label{fig:2d-shapes}Reference domains in two dimensions.}
\end{figure}

Our reference triangle can be seen in \autoref{fig:2d-shapes} and has an
area given by $\int_{-1}^{1}\int_{-1}^{-y}\mathrm{d}x\,\mathrm{d}y = 2$.
A triangle has six symmetries: two rotations, three reflections, and the
identity transformation.  A simple means of realising these symmetries
is to transform into barycentric coordinates
\begin{equation}
  \bm{\lambda} = (\lambda_1,\lambda_2,\lambda_3)^T \quad 0 \leq \lambda_i \leq 1,
  \lambda_1 + \lambda_2 + \lambda_3 = 1,
\end{equation}
which are related to Cartesian coordinates via
\begin{equation}
  \vec{x} = \begin{pmatrix}
    -1 & \hphantom{-}1  & -1\\
    -1 & -1 & \hphantom{-}1\\
  \end{pmatrix}\bm{\lambda},
\end{equation}
where the columns of the matrix can be seen to be the vertices of our
reference triangle.  The utility of barycentric coordinates is that the
symmetric counterparts to a point $\bm{\lambda}$ are given by its unique
permutations.  The number of unique permutations depends on the number
of distinct components of $\bm{\lambda}$ and leads us to the following
three symmetry orbits
\begin{align*}
  \mathcal{S}_1 &= \big(\tfrac{1}{3},\tfrac{1}{3},\tfrac{1}{3}\big), & \abs{\mathcal{S}_1} &= 1,\\
  \mathcal{S}_2(\alpha) &= \Perm(\alpha,\alpha,1-2\alpha), & \abs{\mathcal{S}_2} &= 3,\\
  \mathcal{S}_3(\alpha,\beta) &= \Perm(\alpha,\beta,1-\alpha-\beta), & \abs{\mathcal{S}_3} &= 6,
\end{align*}
where $\alpha$ and $\beta$ are suitably constrained as to ensure the
validity of the resulting coordinates.

It can be easily verified that the orthonormal polynomial basis inside
of our reference triangle is given by
\begin{equation}
  \psi_{ij}(\vec{x}) = \sqrt{2}\hat{P}_i(a)\hat{P}_j^{(2i+1,0)}(b)(1-b)^i,
\end{equation}
where $a = 2(1 + x)/(1 - y) - 1$, and $b = y$ with the objective basis
being given by
\begin{equation}
\mathcal{\tilde{P}}^\phi = \bigl\{\psi_{ij}(\vec{x}) \mid 0 \leq
  i \leq \phi,\; i \leq j \leq \phi - i \bigr\}.
\end{equation}
In the asymptotic limit the cardinality of the objective basis is half
that of the complete basis.  However, the modes of this objective basis
are known not to be completely independent.  Several authors have
investigated the derivation of an optimal quadrature basis on the
triangle.  Details can be found in the papers of Lyness
\cite{lyness1975moderate} and Dunavant \cite{dunavant1985high}.

\paragraph{Quadrilateral.}
Our reference quadrilateral can be seen in \autoref{fig:2d-shapes}.  The
area is simply $\int_{-1}^{1}\int_{-1}^{1}\mathrm{d}x\,\mathrm{d}y = 4$.
A square has eight symmetries: three rotations, four reflections and the
identity transformation.  Applying these symmetries to a point
$(\alpha,\beta)$ with $0 \leq (\alpha,\beta) \leq 1$ will yield a set
$\chi(\alpha,\beta)$ containing its counterparts.  The cardinality of
$\chi$ depends on if any of the symmetries give rise to identical
points.  This can be seen to occur when either $\beta = \alpha$ or
$\beta = 0$.  Enumerating the possible combinations of the above
conditions gives rise to the following four symmetry orbits
\begin{align*}
  \mathcal{S}_1 &= (0, 0), & \abs{\mathcal{S}_1} &=
  1,\\
  \mathcal{S}_2(\alpha) &= \chi(\alpha, 0), & \abs{\mathcal{S}_2} &= 4,\\
  \mathcal{S}_3(\alpha) &= \chi(\alpha, \alpha), & \abs{\mathcal{S}_3} &=
  4,\\
  \mathcal{S}_4(\alpha,\beta) &= \chi(\alpha, \beta), &
  \abs{\mathcal{S}_4} &= 8.
\end{align*}

Trivially, the orthonormal basis inside of our quadrilateral is given by
\begin{equation}
  \psi_{ij}(\vec{x}) = \hat{P}_i(a)\hat{P}_j(b),
\end{equation}
where $a = x$, and $b = y$.  The objective basis is found to be
\begin{equation}
  \mathcal{\tilde{P}}^\phi = \bigl\{\psi_{ij}(\vec{x}) \mid 0 \leq
  i \leq \phi,\; i \leq j \leq \phi - i,\; (i, j) \text{ even} \bigr\},
\end{equation}
with a cardinality one eighth that of the complete basis.

\begin{figure}
  \centering
  \begin{subfigure}[b]{.49\linewidth}
    \centering
    \includegraphics{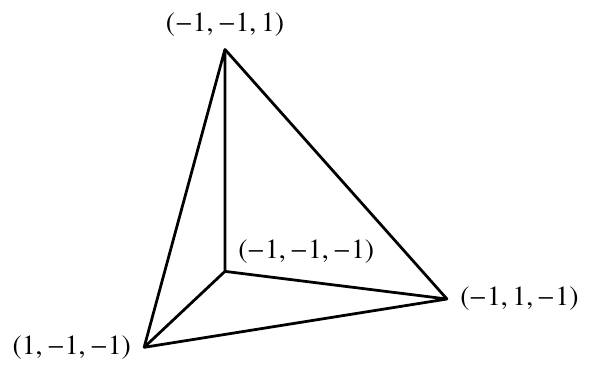}
    \caption{Tetrahedron.}
  \end{subfigure}
  \begin{subfigure}[b]{.49\linewidth}
    \centering
    \includegraphics{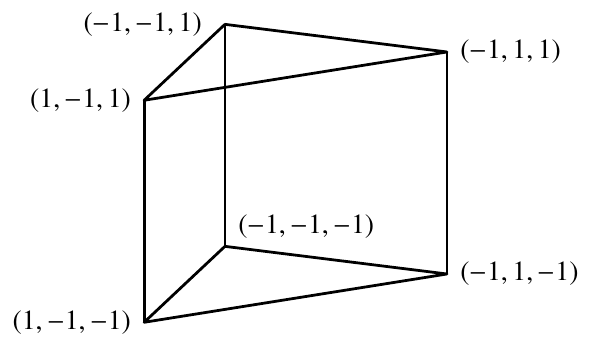}
    \caption{Prism.}
  \end{subfigure}\vskip12pt
  \begin{subfigure}[b]{.49\linewidth}
    \centering
    \includegraphics{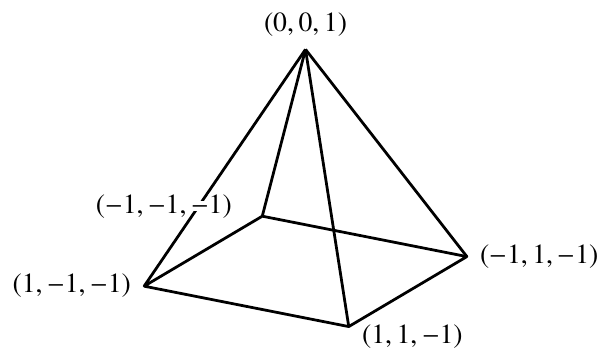}
    \caption{Pyramid.}
  \end{subfigure}
  \begin{subfigure}[b]{.49\linewidth}
    \centering
    \includegraphics{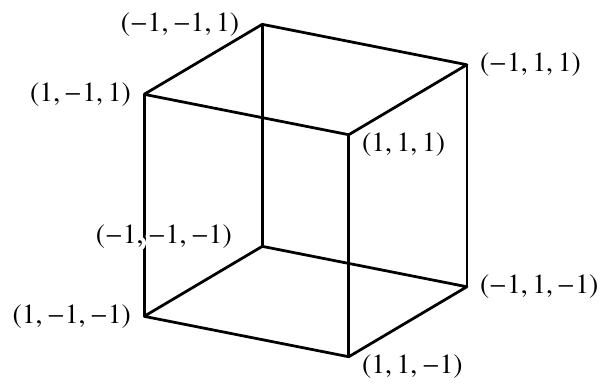}
    \caption{Hexahedron.}
  \end{subfigure}
  \caption{\label{fig:3d-shapes}Reference domains in three dimensions.}
\end{figure}

\paragraph{Tetrahedron.}
Our reference tetrahedron is a right-tetrahedron as depicted in
\autoref{fig:3d-shapes}.  Integrating up the volume we find
$\int_{-1}^{-1}\int_{-1}^{-z}\int_{-1}^{-1-y-z}
\mathrm{d}x\,\mathrm{d}y\,\mathrm{d}z = 4/3$.  A tetrahedron has a total
of 24 symmetries.  Once again it is convenient to work in terms of
barycentric coordinates which are specified for a tetrahedron as
\begin{equation}
  \bm{\lambda} = (\lambda_1,\lambda_2,\lambda_3,\lambda_4)^T \quad 0 \leq \lambda_i \leq 1,
  \lambda_1 + \lambda_2 + \lambda_3 + \lambda_4 = 1,
\end{equation}
and related to Cartesian coordinates via
\begin{equation}
  \vec{x} = \begin{pmatrix}
    -1 & \hphantom{-}1  & -1 & -1\\
    -1 & -1 & \hphantom{-}1 & -1 \\
    -1 & -1 & -1 & \hphantom{-}1
  \end{pmatrix}\bm{\lambda},
\end{equation}
where as with the triangle the columns of the matrix correspond to
vertices of the reference tetrahedron.  Similarly the symmetric
counterparts of $\bm{\lambda}$ are given by its unique permutations.
This leads us to the following five symmetry orbits
\begin{align*}
  \mathcal{S}_1 &=
  \big(\tfrac{1}{4},\tfrac{1}{4},\tfrac{1}{4},\tfrac{1}{4}\big), &
  \abs{\mathcal{S}_1} &= 1,\\
  \mathcal{S}_2(\alpha) &= \Perm(\alpha,\alpha,\alpha,1-3\alpha), &
  \abs{\mathcal{S}_2} &= 4,\\
  \mathcal{S}_3(\alpha) &= \Perm\big(\alpha,\alpha,\tfrac{1}{2} -
  \alpha, \tfrac{1}{2} - \alpha\big), &
  \abs{\mathcal{S}_3} &= 6,\\
  \mathcal{S}_4(\alpha, \beta) &=
  \Perm(\alpha,\alpha,\beta,1-2\alpha-\beta), &
  \abs{\mathcal{S}_4} &= 12,\\
  \mathcal{S}_5(\alpha,\beta,\gamma) &= \Perm(\alpha,\beta,\gamma,
  1-\alpha-\beta-\gamma), & \abs{\mathcal{S}_5} &= 24,
\end{align*}
where $\alpha$, $\beta$, and $\gamma$ are constrained to ensure that $0
\leq \lambda_i \leq 1$ and $\sum_i \lambda_i = 1$.

With some manipulation it can be verified that the orthonormal
polynomial basis inside of our reference tetrahedron is given by
\begin{equation}
  \psi_{ijk}(\vec{x}) = \sqrt{8}\hat{P}_i(a)\hat{P}_j^{(2i+1,0)}(b)
  \hat{P}_k^{(2i + 2j + 2, 0)}(c)(1 - b)^i(1 - c)^{i + j},
\end{equation}
where $a = -2(1 + x)/(y + z) - 1$, $b = 2(1 + y)/(1 - z)$, and $c = z$.
The objective basis is given by
\begin{equation}
\mathcal{\tilde{P}}^\phi = \bigl\{\psi_{ijk}(\vec{x}) \mid 0 \leq
  i \leq \phi, i \leq j \leq \phi - i, j \leq k \leq \phi - i - j \bigr\}.
\end{equation}

\paragraph{Prism.}
Extruding the reference triangle along the $z$-axis gives our reference
prism of \autoref{fig:3d-shapes}.  It follows that the volume is
$\int_{-1}^{1}\int_{-1}^{1}\int_{-1}^{-y}
\mathrm{d}x\,\mathrm{d}y\,\mathrm{d}z = 4$.  There are a total of 12
symmetries.  On account of the extrusion the most natural coordinate
system is a combination of barycentric and Cartesian coordinates:
$(\lambda_1, \lambda_2, \lambda_3, z)$.  Let $\Perm_3$ generate
all of the unique permutations of its first three arguments.  Using this
the six symmetry groups of the prism can be expressed as
\begin{align*}
  \mathcal{S}_1 &= (\tfrac{1}{3},\tfrac{1}{3},\tfrac{1}{3},0), & \abs{\mathcal{S}_1} &= 1,\\
  \mathcal{S}_2(\gamma) &= (\tfrac{1}{3},\tfrac{1}{3},\tfrac{1}{3},\pm\gamma), & \abs{\mathcal{S}_2} &=
  2,\\
  \mathcal{S}_3(\alpha) &= \Perm_3(\alpha,\alpha,1-2\alpha,0), &
  \abs{\mathcal{S}_3} &= 3,\\
  \mathcal{S}_4(\alpha,\gamma) &=
  \Perm_3(\alpha,\alpha,1-2\alpha,\pm\gamma), & \abs{\mathcal{S}_4} &=
  6,\\
  \mathcal{S}_5(\alpha,\beta) &= \Perm_3(\alpha,\beta,1-\alpha-\beta,0),
  & \abs{\mathcal{S}_5} &= 6,\\
  \mathcal{S}_6(\alpha,\beta,\gamma) &=
  \Perm_3(\alpha,\beta,1-\alpha-\beta, \pm\gamma), & \abs{\mathcal{S}_6}
  &= 12,
\end{align*}
where the constraints on $\alpha$ and $\beta$ are identical to those in
a triangle and $0 < \gamma \leq 1$.

Combining the orthonormal polynomial bases for a right-triangle and line
segment yields the orthonormal prism basis
\begin{equation}
  \psi_{ijk}(\vec{x}) = \sqrt{2}\hat{P}_i(a)\hat{P}_j^{(2i + 1, 0)}(b)
  \hat{P}_k(c)(1 - b)^i,
\end{equation}
where $a = 2(1 + x)/(1 - y) - 1$, $b = y$, and $c = z$.  The objective
basis is given by
\begin{equation}
\mathcal{\tilde{P}}^\phi = \bigl\{\psi_{ijk}(\vec{x}) \mid 0 \leq
  i \leq \phi,\; i \leq j \leq \phi - i,\; 0 \leq k \leq \phi - i - j,\;
  k \text{ even}\bigr\}.
\end{equation}

\paragraph{Pyramid.}
Our reference pyramid can be seen in \autoref{fig:3d-shapes} with a
volume determined by
$\int_{-1}^{1}\int_{(z-1)/2}^{(1-z)/2}\int_{(z-1)/2}^{(1-z)/2}
\mathrm{d}x\,\mathrm{d}y\,\mathrm{d}z = 8/3$.  The symmetries are
identical to those of a quadrilateral.  Extending the notation employed
for the quadrilateral we obtain the following symmetry orbits
\begin{align*}
  \mathcal{S}_1(\gamma) &= (0, 0, \gamma), & \abs{\mathcal{S}_1} &=
  1,\\
  \mathcal{S}_2(\alpha,\gamma) &= \chi(\alpha, 0, \gamma), & \abs{\mathcal{S}_2} &= 4,\\
  \mathcal{S}_3(\alpha,\gamma) &= \chi(\alpha, \alpha, \gamma), & \abs{\mathcal{S}_3} &=
  4,\\
  \mathcal{S}_4(\alpha,\beta,\gamma) &= \chi(\alpha, \beta, \gamma), &
  \abs{\mathcal{S}_4} &= 8,
\end{align*}
subject to the constraints that $0 < (\alpha,\beta) \leq (1 - \gamma)/2$
and $-1 \leq \gamma \leq 1$.

Inside of the reference pyramid the orthonormal polynomial basis is
found to be
\begin{equation}
  \psi_{ijk}(\vec{x}) = 2\hat{P}_i(a)\hat{P}_j(b) \hat{P}_k^{(2i + 2j + 2,
    0)}(c)(1 - c)^{i + j},
\end{equation}
where $a = 2x/(1 - z)$, $b = 2y/(1 - z)$, and $c = z$.  The objective
basis is
\begin{equation}
\mathcal{\tilde{P}}^\phi = \bigl\{\psi_{ijk}(\vec{x}) \mid 0 \leq
  i \leq \phi,\; i \leq j \leq \phi - i,\; 0 \leq k \leq \phi - i - j,\;
  (i,j) \text{ even}\bigr\}.
\end{equation}

\paragraph{Hexahedron.}
Our choice of reference hexahedron can be seen in
\autoref{fig:3d-shapes}.  The volume is, trivially,
$\int_{-1}^{1}\int_{-1}^{1}\int_{-1}^{1}
\mathrm{d}x\,\mathrm{d}y\,\mathrm{d}z = 8$.  A hexahedron exhibits
octahedral symmetry with a symmetry number of 48.  The procedure for
determining the orbits similar to that used for the quadrilateral.  We
consider applying these symmetries to a point $(\alpha, \beta, \gamma)$
with $0 \leq (\alpha, \beta, \gamma) \leq 1$ and let the resulting set
of points by given by $\Xi(\alpha, \beta, \gamma)$.  When $\alpha$,
$\beta$, and $\gamma$ are all distinct and greater than zero the set has
a cardinality of 48, as expected.  However, when one or more parameters
are either identical to one another or equal to zero some symmetries
give rise to equivalent points.  This reduces the cardinality of the
set.  Enumerating the various combinations we obtain seven symmetry
orbits
\begin{align*}
  \mathcal{S}_1 &= \Xi(0, 0, 0), & \abs{\mathcal{S}_1} &= 1,\\
  \mathcal{S}_2(\alpha) &= \Xi(\alpha, 0, 0), & \abs{\mathcal{S}_2} &=
  6,\\
  \mathcal{S}_3(\alpha) &= \Xi(\alpha, \alpha, \alpha), &
  \abs{\mathcal{S}_3} &= 8,\\
  \mathcal{S}_4(\alpha) &= \Xi(\alpha, \alpha, 0), & \abs{\mathcal{S}_4}
  &= 12,\\
  \mathcal{S}_5(\alpha, \beta) &= \Xi(\alpha, \beta, 0), &
  \abs{\mathcal{S}_5} &= 24,\\
  \mathcal{S}_6(\alpha, \beta) &= \Xi(\alpha, \alpha, \beta), &
  \abs{\mathcal{S}_6} &= 24,\\
  \mathcal{S}_7(\alpha, \beta, \gamma) &= \Xi(\alpha, \beta, \gamma), &
  \abs{\mathcal{S}_7} &= 48.
\end{align*}
Trivially, the orthonormal basis inside of our reference hexahedron is
given by
\begin{equation}
  \psi_{ijk}(\vec{x}) = \hat{P}_i(a)\hat{P}_j(b)\hat{P}_k(c),
\end{equation}
where $a = x$, $b = y$, and $c = z$.  The objective basis is
\begin{equation}
\mathcal{\tilde{P}}^\phi = \bigl\{\psi_{ijk}(\vec{x}) \mid 0 \leq
  i \leq \phi,\; i \leq j \leq \phi - i,\; j \leq k \leq \phi - i - j,\; (i,
  j, k) \text{ even} \bigr\}.
\end{equation}

\section{Methodology}
\label{sec:meth}

Our methodology for identifying symmetric quadrature rules is a
refinement of that described by Witherden and Vincent
\cite{witherden2013analysis} for triangles.  This method is, in turn, a
refinement of that of Zhang et al. \cite{zhang2009set}.

To derive a quadrature rule four input parameters are required: the
reference domain, $\vec{\Omega}$, the number of quadrature points $N_p$,
the target rule strength, $\phi$, and a desired runtime, $t$.  The algorithm
begins by computing all of the possible symmetric decompositions of
$N_p$.  The result is a set of vectors satisfying the relation
\begin{equation}
  N_p = \sum^{N_s}_{j = 1} n_{ij}|\mathcal{S}_j|,
\end{equation}
where $N_s$ is the number of symmetry orbits associated with the domain
$\vec{\Omega}$, and $n_{ij}$ is the number of orbits of type $j$ in the
$i$th decomposition.  Finding these involves solving the constrained
linear Diophantine equation outlined in \autoref{sec:shapes}.  It is
possible for this equation to have no solutions.  As an example we
consider the case when $N_p = 44$ for a triangular domain.  From the
symmetry orbits we have
\[
  N_p = n_1|\mathcal{S}_1| + n_2|\mathcal{S}_2| +
  n_3|\mathcal{S}_3| = n_1 + 3n_2 + 6n_3,
\]
subject to the constraint that $n_1 \in \{0, 1\}$.  This restricts $N_p$
to be either a multiple of three or one greater.  Since forty-four is
neither of these we find the equation to have no solutions.  Therefore,
we conclude that there can be no symmetric quadrature rules inside of a
triangle with forty-four points.

Given a decomposition we are interested in finding a set of orbital
parameters and weights that minimise the error associated with
integrating the objective basis on $\vec{\Omega}$.  This is an example
of a nonlinear least squares problem.  A suitable method for solving
such problems is the Levenberg-Marquardt algorithm (LMA).  The LMA is an
iterative procedure for finding a set of parameters that correspond to a
local minima of a set of functions.  The minimisation process is not
always successful and is dependent on an initial guess of the
parameters.  Within the context of quadrature rule derivation
minimisation can be regarded as successful if $\xi(\phi) \sim \epsilon$
where $\epsilon$ represents machine precision.

Let us denote the number of parameters associated with symmetry orbit
$\mathcal{S}_i$ as $\big\llbracket\mathcal{S}_i\big\rrbracket$.  Using
this we can express the total number of degrees of freedom associated
with decomposition $i$ as
\begin{equation}
  \sum^{N_s}_{j = 1} \big\{n_{ij}\big\llbracket\mathcal{S}_i\big\rrbracket + n_{ij}\big\},
\end{equation}
with the second term accounting for the presence of one quadrature
weight associated with each symmetry orbit.  From the list of orbits
given in \autoref{sec:shapes} we expect the weights contribute
approximately one third of the degrees of freedom.  This is not an
insignificant fraction.  One way of eliminating the weights is to treat
them as dependent variables.  When the points are prescribed the right
hand side of \autoref{eq:quad} becomes linear with respect to the
unknowns---the weights.  In general, however, the number of weights will
be different from the number of polynomials in the objective basis.  It
is therefore necessary to obtain a least squares solution to the systen.
Linear least squares problems can be solved directly through a variety
of techniques.  Perhaps the most robust numerical scheme is that of
singular value decomposition (SVD).  Thus, at the cost of solving a
small linear least squares problem at each LMA iteration we are able to
reduce the number of free parameters to
\begin{equation}
  \sum^{N_s}_{j = 1} n_{ij}\big\llbracket\mathcal{S}_i\big\rrbracket.
\end{equation}
Such a modification has been found to greatly reduce the number of
iterations required for convergence.  This reduction more than offsets
the marginally greater computational cost associated with each
iteration.

Previous works \cite{zhang2009set, kubatko2013pri, taylor2005several}
have emphasised the importance of picking a `good' initial guess to seed
the LMA.  To this end several methodologies for seeding orbital
parameters have been proposed.  The degree of complexity associated with
such strategies is not insignificant.  Further, it is necessary to
device a separate strategy for each symmetry orbit.  Our experience,
however, suggests that the choice of decomposition is far more important
than the initial guess in determining whether minimisation will be
successful.  For larger values of $N_p$ we note that many
decompositions---especially those for prisms and pyramids---are
pathological.  As an example of this we consider searching for an $N_p =
80$ point rule inside of a prism where there are $2\,380$ distinct
symmetrical decompositions.  One such decomposition is $N_p =
40\abs{\mathcal{S}_2}$ where all points lie in a single column down the
middle of the prism.  Since there is no variation in either $x$ or $y$
it is not possible to obtain a rule of strength $\phi \geq 1$.  Hence,
the decomposition can be dismissed without further consideration.

A presentation of our method in pseudocode can be seen in
\autoref{alg:meth}.  When the objective basis functions in
$\mathcal{\tilde{P}}^\phi$ are orthonormal \autoref{eq:obint} states that
the integrand of all non-constant modes is zero.  We can exploit this to
simplify the computation of $b_i$.  The purpose of \textsc{ClampOrbit}
is to enforce the constraints associated with a given orbit to ensure
that all points remain inside of the domain.

\begin{algorithm}
  \caption{\label{alg:meth}Procedure for generating symmetric quadrature
  rules of strength $\phi$ with $N_p$ points inside of a domain.}
  \begin{algorithmic}[1]
    \Procedure{FindRules}{$N_p,\phi,t$}
      \ForAll{decompositions of $N_p$}
        \State $t_0 \gets \Call{CurrentTime}$
        \Repeat
        \State $\mathcal{R} \gets \Call{SeedOrbits}$
        \Comment{Initial guess of points}
        \State $\xi \gets \Call{LMA}{\textsc{RuleResid}, \mathcal{R}}$
        \If{$\xi \sim \epsilon$}
          \Comment{If minimisation was successful}
          \State \textbf{save} $\mathcal{R}$
          \EndIf
        \Until{$\Call{CurrentTime}{}() - t_0 > t$}
      \EndFor
    \EndProcedure
    \Statex
    \Function{RuleResid}{$\mathcal{R}$}
      \ForAll{$p_i \in \mathcal{\tilde{P}}^\phi$}
      \Comment{For each basis function}
        \State $b_i \gets \int_{\vec{\Omega}} p_i(\vec{x})\,
        \mathrm{d}\vec{x}$
        \ForAll{$r_j \in \mathcal{R}$}
        \Comment{For each orbit}
          \State $r_j \gets \Call{ClampOrbit}{r_j}$
          \Comment{Ensure orbital parameters are valid}
          \State $A_{ij} \gets 0$
          \ForAll{$\vec{x}_k \in \Call{ExpandOrbit}{r_j}$}
            \State $A_{ij} \gets A_{ij} + p_i(\vec{x}_k)$
          \EndFor
        \EndFor
      \EndFor
      \State $\omega \gets b / A$
      \Comment{Use SVD to determine the weights}
      \State \textbf{return} $A\omega - b$
      \Comment{Compute the residual}
    \EndFunction
  \end{algorithmic}
\end{algorithm}

\section{Implementation}
\label{sec:impl}

We have implemented the algorithms outlined above in a C++11 program
called \emph{polyquad}.  The program is built on top of the Eigen
template library \cite{eigenweb} and is parallelised using MPI.  It is
capable of searching for quadrature rules on triangles, quadrilaterals,
tetrahedra, prisms, pyramids, and hexahedra.  All rules are guaranteed
to be symmetric having all points inside of the domain.  Polyquad can
also, optionally, filter out rules possessing negative weights.
Further, functionality exists, courtesy of MPFR \cite{fousse2007mpfr},
for refining rules to an arbitrary degree of numerical precision and for
evaluating the truncation error of a ruleset.

The source code for polyquad is available under the terms of the GNU
General Public License v3.0 and can be downloaded from
\url{https://github.com/vincentlab/Polyquad}.

\section{Rules}
\label{sec:rules}

Using polyquad we have derived a set of quadrature rules for each of the
reference domains in \autoref{sec:shapes}.  All rules are completely
symmetric, possess only positive weights, and have all points
inside of the domain.  It is customary in the literature to
refer to quadratures with the last two attributes as being ``PI'' rules.  As
polyquad attempts to find an ensemble of rules it is necessary to have a
means of differentiating between otherwise equivalent formulae.  In
constructing this collection the truncation term $\xi(\phi + 1)$ was
employed with the rule possessing the smallest such term being chosen.
The number of points $N_p$ required for a rule of strength $\phi$ can be
seen in \autoref{tab:rules}.  The rules themselves are provided as
electronic supplementary material and have been refined to 38 decimal
places.

\begin{table}
  \centering
  \caption{\label{tab:rules}Number of points $N_p$ required for a fully
    symmetric quadrature rule with positive weights of strength $\phi$
    inside of the six reference domains.  Rules with underlines
    represent improvements over those found in the literature (see text).}
  \begin{tabular}{rrrrrrr} \toprule
    & \multicolumn{6}{c}{$N_p$} \\ \cmidrule{2-7}
    $\phi$ & Tri & Quad & Tet & Pri & Pyr & Hex\\ \midrule
    1  & 1  & 1  & 1  & 1  & 1  & 1  \\
    2  & 3  & 4  & 4  & 5  & 5  & 6  \\
    3  & 6  & 4  & 8  & 8  & 6  & 6  \\
    4  & 6  & 8  & 14 & 11 & 10 & 14 \\
    5  & 7  & 8  & 14 & 16 & \underline{15} & 14 \\
    6  & 12 & 12 & 24 & \underline{28} & \underline{24}  & \underline{34} \\
    7  & 15 & 12 & \underline{35} & \underline{35} & \underline{31} &
    \underline{34} \\
    8  & 16 & \underline{20} & 46 & \underline{46} & \underline{47} & \underline{58} \\
    9  & 19 & \underline{20} & \underline{59} & \underline{60} & \underline{62} &
    \underline{58} \\
    10 & 25 & \underline{28} & 81 & \underline{85} & \underline{83} &
    \underline{90} \\
    11 & 28 & \underline{28} &   &   &   & \\
    12 & 33 & \underline{37} &   &   &   & \\
    13 & 37 & \underline{37} &   &   &   & \\
    14 & 42 & \underline{48} &   &   &   & \\
    15 & 49 & \underline{48} &   &   &   & \\
    16 & 55 & 60 &   &   &   & \\
    17 & 60 & 60 &   &   &   & \\
    18 & 67 & \underline{72} &   &   &   & \\
    19 & 73 & \underline{72} &   &   &   & \\
    20 & 79 & \underline{85} &   &   &   & \\
    \bottomrule
  \end{tabular}
\end{table}

From the table we note that several of the rules appear to improve over
those in the literature.  We consider a rule to be an improvement when
it either requires fewer points than any symmetric rule described in
literature or when existing symmetric rules of this strength are not PI.
We note that many of the rules presented by Dunavant for quadrilaterals
\cite{dunavant1985economical} and hexahedra \cite{dunavant1986efficient}
possess either negative weights or have points outside of the domain.
Using polyquad in quadrilaterals we were able to identify PI rules with
point counts less than or equal to those of Dunavant at strengths $\phi
= 8,9,18,19,20$.  In tetrahedra we were able to reduce the number of
points required for $\phi = 7$ and $\phi = 9$ by one and two,
respectively, compared with Zhang et al. \cite{zhang2009set}.
Furthermore, in prisms and pyramids rules requiring significantly fewer
points than those in literature were identified.  As an example, the
$\phi = 9$ rule of \cite{kubatko2013pri} inside of a prism requires 71
points compared with just 60 for the rule identified by polyquad.
Additionally, both of the $\phi = 10$ rules for prisms and pyramids
appear to be new.

\section{Conclusions}
\label{sec:conclusions}

We have presented a methodology for identifying symmetric quadrature rules
on a variety of domains in two and three dimensions. Our scheme does not
require any manual intervention and is not restricted to any particular
topological configuration inside of a domain. Additionally, it is also
capable of generating an ensemble of rules. We have further provided an
open source implementation of our method in C++11, and used it to generate
a complete set of symmetric quadrature rules that are suitable for use in
finite element solvers. All rules possess purely positive weights and have
all points inside the domain. Many of the rules appear to be new, and an
improvement over those tabulated in the literature.

\section*{Acknowledgements}
The authors would like to thank the Engineering and Physical Sciences
Research Council for their support via a Doctoral Training Grant and
an Early Career Fellowship (EP/K027379/1).

\end{document}